\documentclass[12pt]{article}
\usepackage{amsfonts,amssymb,amsmath}

\newcommand\CC{\mathbb C}
\newcommand\beq{\begin{equation}}
\newcommand\eeq{\end{equation}}
\newtheorem{theorem}{Theorem}
\newtheorem{prop}{Proposition}
\newtheorem{remark}{Remark}
\begin{document}
\title{Generalized analytic functions, Moutard-type transforms and holomorphic 
maps
\thanks{The main part of the work was fulfilled during the visit of the 
first author to the IHES, France in November 2015. The first author was partially supported 
by the Russian Foundation for Basic Research, grant 13-01-12469 ofi-m2 and by the program 
``Fundamental problems of nonlinear dynamics'', RAS.}}
\author{P.G. Grinevich
\thanks{Landau Institute of Theoretical Physics, Chernogolovka, 142432,
Russia; Lomonosov Moscow State University, Moscow, 119991, Russia; e-mail: pgg@landau.ac.ru} \and R.G. Novikov\thanks
{CNRS (UMR 7641), Centre de Math\'ematiques Appliqu\'ees, 
\'Ecole Polytechnique, 91128, Palaiseau, France;
IEPT RAS, 117997, Moscow, Russia;
e-mail: novikov@cmap.polytechnique.fr}}
\date{}
\maketitle
\begin{abstract}
We continue the studies of Moutard-type transform for generalized analytic functions 
started in \cite{GN}. In particular, we suggest an interpretation of generalized analytic 
functions as spinor fields and show that in the framework of this approach Moutard-type
transforms for aforementioned functions commute with holomorphic changes of variables.
\end{abstract}

Key words: generalized analytic functions, spinors, Moutard transforms.

\medskip

We study the basic pair of conjugate equations of the generalized analytic 
function theory:  
\begin{align}
\label{g1}
&\partial_{\bar z} \psi = u \bar \psi \ \ \mbox{in} \ \ D,\\
\label{g2}
&\partial_{\bar z} \psi^+ = -\bar u \bar \psi^+ \ \ \mbox{in} \ \ D,
\end{align}
where $D$ is an open domain in $\CC$,  $u=u(z)$ is a given function in $D$, 
$\partial_{\bar z}=\partial/\partial\bar z$; see \cite{V}. Here and below, 
the notation $f=f(z)$ does not mean that $f$ is holomorphic. 

A new progress in the theory of generalized analytic functions was obtained very recently 
in \cite{GN} by showing that Moutard-type transforms can be applied to the pair 
of equations  (\ref{g1}), (\ref{g2}). Note that ideas of Moutard-type transforms 
were developed and successfully used in the differential geometry, in the soliton theory 
in dimension 2+1, and in the spectral theory in dimension 2, see  \cite{GN} for further 
references. In particular, our work \cite{GN} was essentially stimulated by recent articles by 
I.A. Taimanov \cite{T1}, \cite{T2} on the Moutard-type transforms for the Dirac operators
in the framework of the soliton theory in the dimension $2+1$. On the other hand, we were 
strongly motivated by some open problems of two-dimensional inverse scattering at fixed energy, 
where equation  (\ref{g1}) arises as the $\bar\partial$-equation in spectral parameter.

A simple Moutard-type transform $\mathcal{M}=\mathcal{M}_{u,f,f^+}$ for the pair 
of conjugate equations (\ref{g1}), (\ref{g2}) is given by the formulas
(see \cite{GN}):
\beq
\label{m3}
\tilde u = \mathcal{M} u= u + \frac{f\overline{f^+}}{\omega_{f,f^+}},
\eeq
\beq
\label{m1}
\tilde\psi=\mathcal{M} \psi=
\psi-\frac{\omega_{_{\psi,f^+}}}{\omega_{f,f^+}}\,f , \ \ 
\tilde{\psi}^+=  \mathcal{M} \psi^+=  \psi^+ - \frac{\omega_{f,\psi^+}}{\omega_{f,f^+}}\, f^+,
\eeq
where $f$ and $f^+$ are some fixed solutions of equations 
(\ref{g1}) and (\ref{g2}), respectively, $\psi$ and  
$\psi^+$ are arbitrary solutions of (\ref{g1}) and (\ref{g2}),
respectively, and  $\omega_{\psi,\psi^+}=\omega_{\psi,\psi^+}(z)$ denotes 
imaginary-valued function defined by:
\beq
\label{k1}
\partial_{z} \omega_{\psi,\psi^+} =\psi\psi^+, \ \ 
\partial_{\bar z} \omega_{\psi,\psi^+} =-\overline{\psi\psi^+} \ \ \mbox{in} \ \ D,
\eeq
where this definition is self-consistent, at least, for simply connected $D$, 
whereas a pure imaginary integration constant may depend on concrete 
situation. The point is that the functions $\tilde\psi$, $\tilde\psi^+$ defined in (\ref{m1})
satisfy the conjugate pair of Moutard-transformed equations:
\begin{align}
\label{g3}
&\partial_{\bar z} \tilde\psi= \tilde u\, \overline{\tilde\psi}& \ \ &\mbox{in} \ \ D, \\ 
\label{g4}
&\partial_{\bar z} \tilde\psi^+= -\overline{\tilde u}\, \overline{\tilde\psi^+}& \ \ 
&\mbox{in} \ \ D.
\end{align}
In addition, we have also the following new important result: 
\begin{prop}
\label{th:1}
For a simple Moutard transform (\ref{m3}),(\ref{m1}) the following formula holds:
\beq
\label{eq:trans:pot}
\omega_{\tilde\psi,\tilde\psi^+}=\frac{\omega_{\psi,\psi^+}\omega_{f,f^+}-\omega_{\psi,f^+}
\omega_{f,\psi^+}}{\omega_{f,f^+}}+c_{\tilde\psi,\tilde\psi^+},
\eeq 
where $c_{\tilde\psi,\tilde\psi^+}$ is an imaginary constant.
\end{prop}
In order to apply the Moutard-type transforms to studies of generalized-analytic functions with 
contour poles, we need to study, in particular, compositions of the former transforms and 
holomorphic maps. 

Consider a holomorphic bijection $W$:
\beq
\label{c1}
W:D\rightarrow D_*, \ \ z\rightarrow \zeta(z),
\eeq
$$
W^{-1}:D_*\rightarrow D, \ \ \zeta\rightarrow z(\zeta),
$$
where $D$ is the domain in (\ref{g1}), (\ref{g2}).

If we treat $\psi(z)$, $\psi^+(z)$ as scalar fields  in equations (\ref{g1}), (\ref{g2}), 
then the conjugate property of these equations is not invariant with respect to holomorphic 
bijections . In the next theorem we give the proper transformation formulas 
for the conjugate pair of equations  (\ref{g1}), (\ref{g2}) with respect to holomorphic 
bijections:

\begin{theorem}
\label{th:2} Let $W$ be a holomorphic bijection as in (\ref{c1}). Let 
\begin{align}
\label{c2}
&u_*(\zeta)=u(z(\zeta))\,\sqrt{\frac{\partial z}{\partial \zeta}\,
\frac{\partial\bar z}{\partial\bar\zeta} }= u(z(\zeta))\, \left|\frac{\partial z}{\partial \zeta}
\right|,  \\ 
&\psi_*(\zeta)=\psi(z(\zeta))\,\sqrt{\frac{\partial z}{\partial \zeta}},\ \ \label{c3}
\psi^+_*(\zeta)=\psi^+(z(\zeta))\,\sqrt{\frac{\partial z}{\partial \zeta}}, 
\end{align}
where $u(z)$, $\psi(z)$,  $\psi^+(z)$ are the same that in equations (\ref{g1}), 
(\ref{g2}). Then:
\begin{align}
&\partial_{\bar\zeta} \psi_* = u_* \bar \psi_* \ \ \mbox{in} \ \ D_*, \label{c4}\\
&\partial_{\bar\zeta} \psi^+_* = -\bar u_* \bar \psi^+_* \ \ \mbox{in} \ \ D_*. \label{c5}
\end{align}
In addition, 
\beq
\label{c6}
\omega_{\psi_*,\psi_*^+}(\zeta)=\omega_{\psi,\psi^+}(z(\zeta)),
\eeq
where $\omega$ is defined according to (\ref{k1}).
\end{theorem}
\begin{remark}
Formulas (\ref{c2}),  (\ref{c3}) have the following natural interpretation:
$\psi(z)$,  $\psi^+(z)$ can be treated as spinors, i.e. differential forms of the type 
$\left(\frac{1}{2},0\right)$, and $u$ can be treated as differential form of the type 
$\left(\frac{1}{2},\frac{1}{2}\right)$. The corresponding forms can be written as:
\beq
u=u(z)\sqrt{dzd\bar z}, \ \ \psi=\psi(z)\sqrt{dz}, \ \  \psi^+=\psi^+(z)\sqrt{dz}.
\eeq
It is very natural because the generalized analytic function equation (\ref{g1}) 
can be viewed as a special reduction of the two-dimensional Dirac system, see, for example \cite{GN}.
\end{remark}
Theorem \ref{th:2} implies that $W$ in (\ref{c1}) generates a map of the conjugate pair 
of equations (\ref{g1}), (\ref{g2}) into the conjugate pair of equations (\ref{c4}), 
(\ref{c5}). We also denote the latter map by $W$. Using this interpretation of $W$ we obtain 
the following result:
\begin{theorem}
\label{th:3}
The following formula holds: 
\beq
\mathcal{M}_{u_*,f_*,f^+_*} \circ W = W\circ \mathcal{M}_{u,f,f^+},
\eeq
where $\mathcal{M}_{u,f,f^+}$ and $\mathcal{M}_{u_*,f_*,f^+_*}$ are defined according to formulas  
(\ref{m3}), (\ref{m1}), and $u_*$, $f_*$, $f_*^+$, $\omega_{\psi_*,\psi_*^+}$ 
are defined according to (\ref{c2}),  (\ref{c3}),  (\ref{c6}). 
\end{theorem}
Proposition~\ref{th:1} and Theorems~\ref{th:2} and \ref{th:3} can be proved by direct calculations.

In the framework of the Moutard transform approach, using Theorem~\ref{th:3} we reduce local 
studies of generalized analytic functions with contour pole at a real-analytic curve to the 
case of contour pole at a straight line. These studies will be continued in a subsequent paper.

\end{document}